\documentclass[12pt,draftcls,onecolumn]{IEEEtran}
\usepackage{fancyhdr}

\usepackage{amsfonts}
\usepackage{mathrsfs}

\usepackage{colortbl}
\usepackage{amssymb,amsmath,graphicx}

\usepackage{caption}
\usepackage{subfigure}
\usepackage{float}
\usepackage[outdir=./]{epstopdf}
\graphicspath{{../eps/}}
\DeclareGraphicsExtensions{.eps}

\newtheorem{Remark}{\bf Remark}[section]

\newtheorem{Definition}{\bf Definition}[section]

\newenvironment{Proof}{\noindent{\em Proof:\/}}{\hfill $\Box$\par}

\newtheorem{Theorem}{\bf Theorem}[section]

\newtheorem{Lemma}{\bf Lemma}[section]

\newtheorem{Property}{Property}[section]
\newtheorem{Assumption}{\bf Assumption}[section]

\newcommand{\R}{{\mathbb R}}

\newcommand{\EQQ}{\begin{eqnarray*}}
\newcommand{\ENN}{\end{eqnarray*}}
\newcommand{\EQ}{\begin{eqnarray}}
\newcommand{\EN}{\end{eqnarray}}
\newcommand{\bd}{\begin{Definition}
\newcommand{\mbox{col}}{\mbox{col}}
\begin{rm} }
\newcommand{\ed}{ \end{rm}
\end{Definition} }

\begin{document}

\title{\LARGE \bf
Cooperative Robust Output Regulation Problem for Discrete-Time Linear Time-Delay Multi-Agent Systems
}

\author{Yamin~Yan~and~Jie~Huang
\thanks{This work has been supported by the Research Grants Council of the Hong Kong Special Administration
Region under grant No. 14219516.}
\thanks{Yamin Yan and Jie Huang  are with Department of Mechanical and Automation
Engineering, The Chinese University of Hong Kong, Shatin, N.T., Hong
Kong. E-mail: ymyan@mae.cuhk.edu.hk, jhuang@mae.cuhk.edu.hk}
}

\maketitle

\thispagestyle{fancy}
\fancyhead{}
\lhead{}
\lfoot{}
\cfoot{}
\rfoot{}
\renewcommand{\headrulewidth}{0pt}
\renewcommand{\footrulewidth}{0pt}

\begin{abstract}
In this paper, we study the cooperative robust output
regulation problem for discrete-time linear multi-agent systems
with both communication and input delays by distributed
internal model approach. We first introduce the distributed internal model
for discrete-time multi-agent systems with both
communication and input delays. Then, we define so-called auxiliary system and auxiliary augmented system.
Finally, we solve our problem by showing, under some standard assumptions,  that if a distributed state
feedback control or a  distributed output feedback control solves the robust output regulation problem of the
auxiliary system, then the same control law  solves the cooperative robust output regulation problem of the
original multi-agent systems.
\end{abstract}

\section{INTRODUCTION}

The cooperative output regulation problem aims to design a control law for a multi-agent system to drive the tracking error of each follower to the origin asymptotically while rejecting a class of external disturbances.  The problem is interesting because its formulation includes the leader-following consensus, synchronization or formation as special cases. Like the output regulation problem of a single linear system \cite{Da,Fr,FW}, there are two approaches to handling the cooperative output regulation problem of multi-agent systems. The first one is called feedforward design
\cite{Su1}. This approach makes use of the solution of the regulator
equations and a distributed observer to design an appropriate feedforward term to exactly cancel
the steady-state  tracking error. The second one is called
distributed internal model design \cite{Su3,Jiang1}. This approach employs
a distributed internal model to convert the cooperative output regulation problem of an uncertain multi-agent system to a
simultaneous eigenvalue assignment problem of a multiple augmented system composed of the
given multi-agent system and the distributed internal model. The internal
model approach has at least two advantages over the feedforward design approach in that it can tolerate
perturbations of the plant parameters, and it does not need to solve the regulator equations.

Recently,
the study on the cooperative output regulation problem has been extended to linear multi-agent systems with time-delay and/or communication delay.
Specifically,  the cooperative output regulation problem for linear continuous-time multi-agent systems with time-delay was studied
in \cite{CCC14} via the distributed observer approach and in \cite{Lu2015} via the distributed internal model approach.
The cooperative output regulation problem for linear discrete-time multi-agent systems with time-delay was studied
in \cite{IET} via the discrete distributed observer approach. Since
the discrete distributed observer approach cannot handle the model uncertainties and the control law has to rely on the solution to the discrete regulator equations, we will  further develop a distributed internal model approach to deal with the cooperative output regulation problem for discrete-time multi-agent systems with both input  and communication delays.

To solve our problem, we will first introduce a distributed internal model for linear discrete-time time-delay  multi-agent systems.
This distributed internal model together with the given multi-agent system defines a so-called auxiliary augmented system.
We will show that, if the communication network of the multi-agent system is connected, then our original problem can be converted to the stabilization problem
of the auxiliary augmented system. Due to this result, it suffices to stabilize the auxiliary augmented system via either static state feedback control law or dynamic output feedback control law.
It is noted that the stabilization problem of the auxiliary augmented system is challenging for two reasons. First, the auxiliary augmented system
is also a time-delay system, and, second, it is also subject to communication constraints. We have managed to overcome these challenges by both distributed
static state feedback control law and distributed dynamic output feedback control law.

Technically, our approach is related to the references \cite{Lu2015} and \cite{JSSC}.
The reference \cite{Lu2015}   deals  with the robust output regulation problem for continuous-time linear multi-agent systems with both communication and input delays. Our current work can be viewed as
a discrete-time analog of the framework in \cite{Lu2015}. Since, for time-delay systems, the regulator equations for continuous-time systems and discrete-time systems
are somehow different and the stabilization techniques for continuous-time systems and discrete-time systems are also different, an independent study on the
 discrete-time multi-agent delay systems is necessary. On the other hand, the reference \cite{JSSC} studies the robust output regulation problem for a single linear system with both communication and input delays. The current paper can be viewed as an extension of the results of \cite{JSSC} to multi-agent systems. Compared with \cite{JSSC}, the main technical challenge is to find a distributed control law to satisfy the communication constraints. It is also noted that
the cooperative robust output regulation problem for  delay-free uncertain discrete-time multi-agent systems was also considered by distributed state feedback control in \cite{Liang2014}.
However, for delay-free systems, the regulator equations for  continuous-time systems and discrete-time systems are the same, and the stabilization techniques for
the auxiliary augmented system of the continuous-time systems and discrete-time systems are also the same. The techniques used in \cite{Su3} for continuous-time multi-agent systems
is directly applicable to the discrete-time multi-agent systems.
It is also worth mentioning some other references relevant to the topic of this paper can be found in \cite{chen1}, \cite{chen2}, and \cite{you2013}.

The rest of this paper is organized as follows.
Section \ref{preliminary} formulates our problem. Section \ref{frame1} defines the distributed internal model and the auxiliary augmented system and
presents a framework for converting our original problem to the stabilization problem
of the auxiliary augmented system. Section \ref{result} establishes the main result. An example is used to illustrate our design in Section \ref{example}. Finally the paper is closed with some concluding remarks in Section \ref{concl}.

{\bf Notation.} $\sigma(A)$ denotes the spectrum of a square matrix $A$. For $X_i \in \R^{n_i}$, $i = 1, \dots, m$,
$\mbox{col} (X_1, \dots, X_m) = [X^T_1, \dots, X^T_m]^T$. For $X = [X_1, \cdots, X_m]$ where $X_i \in \R^{n \times 1}$,
$\mbox{vec}(X) = \mbox{col} (X_1, \dots, X_m)$.
For some nonnegative integer $r$, $I[-r,0]$ denotes the set of integers $\{-r,-r+1,\cdots,0\}$ and $\mathcal{C}(I[-r,0],\R^{n})$ denotes the set of functions mapping the integer set $I[-r,0]$ into $\R^{n}$. $\mathbb{Z}^+ =  \{0, 1, \cdots \}$. $\otimes$ denotes the Kronecker product of matrices.

\section{PROBLEM FORMULATION AND
PRELIMINARIES}\label{preliminary}
\subsection{Graph}

A digraph $\mathcal{G}=(\mathcal{V},\mathcal{E})$ consists of a node set $\mathcal{V}=\{1,\cdots,N\}$ and an edge set $\mathcal{E}\subseteq \mathcal{V}\times \mathcal{V}$. An edge of $\mathcal{E}$ from node $i$ to node $j$ is denoted by $(i,j)$, where the nodes $i$ and $j$ are called the parent node and the child node of each other, and the node $i$ is also called a neighbor of the node $j$. Let $\mathcal{N}_i=\{j,(j,i) \in \mathcal{E}\}$ denote the subset of $\mathcal{V}$ which consists of all the neighbors of the node $i$. Edge $(i,j)$ is called undirected if $(i,j)\in \mathcal{E}$ implies that $(j,i)\in \mathcal{E}$. The graph is called undirected if every edge in $\mathcal{E}$ is undirected. If there exists a set of edges $\{(i_1,i_2),\cdots,(i_k,i_{k+1})\}$ in the digraph $\mathcal{G}$, then $i_{k+1}$ is said to be reachable from node $i_1$. A digraph $\mathcal{G}_s=(\mathcal{V}_s,\mathcal{E}_s)$, where $\mathcal{V}_s\subseteq \mathcal{V}$ and $\mathcal{E}_s \subseteq \mathcal{E}\cap (\mathcal{V}_s \times \mathcal{V}_s)$, is a subgraph of the digraph $\mathcal{G}=(\mathcal{V},\mathcal{E})$.
A weighted adjacency matrix of $\mathcal{G}$ is a square matrix denoted by $\mathcal{A}=\left[a_{ij}\right] \in R^{N\times N}$ such that, for $i,j = 1, \cdots, N$,  $a_{ii}=0$, $a_{ij}>0$ $\Leftrightarrow$ $(j,i)\in \mathcal{E}$, and $a_{ij} = a_{ji}$ if $(i,j)$  is undirected.
The Laplacian matrix of a digraph $\mathcal{G}$ is denoted by $\mathcal{L}=[l_{ij}]\in R^{N\times N}$, where $l_{ii}=\sum_{j=1}^Na_{ij}$, and $l_{ij}=-a_{ij}$ if $i\neq j$. More detailed exposition on graph theory can be found in \cite{Graph}.
\vspace{-6pt}

\subsection{Problem Formulation}

In this paper, we consider the cooperative robust output regulation problem for discrete-time linear uncertain time-delay systems of the following form:
\begin{equation}\label{system0}
\begin{split}
x_i(t+1)=&(A+\delta_i A)x_i(t)+(B+\delta_i B)u_i(t-r_{con})+(E+\delta_i E)v(t),\ t\in \mathbb{Z}^+,\\
y_i(t)=&(C+\delta_i C)x_i(t),\ t\in \mathbb{Z}^+,\\
\end{split}
\end{equation}
where $x_i(t)\in \R^{n_i}$, $y_i(t)\in \R^{p_i}$, and $u_i(t)\in \R^{m_i}$ are the system state, measurement output, and control input of the $i^{th}$ subsystem, $r_{con}\in \mathbb{Z}^+$ is the input delay, and $v(t)\in \R^q$ is the exogenous signal representing the reference input to be tracked or/and disturbance to be rejected and is assumed to be generated by the exosystem of the form
\begin{equation}\label{exo}
v(t+1)=Sv(t), \ t\in\mathbb{Z}^+,
\end{equation}
where $S\in \R^{q\times q}$ is a constant matrix.

The regulated output for each subsystem is defined as
\begin{equation}
e_i(t)=y_i(t)-y_0(t), \ i=1,\cdots, N,
\end{equation}
where $y_0(t)=-Fv(t)$.

In \eqref{system0}, the matrices $A$, $B$, $E$, $C$ represent the nominal part of $i^{th}$ plant, while the matrices $\delta_i A$, $\delta_i B$, $\delta_i E$ and $\delta_i C$, represent the uncertain part of the $i^{th}$ plant. For convenience, we denote the system uncertain parameters by a vector
\begin{equation*}
w = \left[
                                                                                                         \begin{array}{c}
                                                                                                           \text{vec}(\delta_1 A,\dots,\delta_N A)\\
                                                                                                           \text{vec}(\delta_1 B,\dots,\delta_N B)\\
                                                                                                           \text{vec}(\delta_1 E,\dots,\delta_N E)\\
                                                                                                           \text{vec}(\delta_1 C,\dots,\delta_N C)\\
                                                                                                         \end{array}
                                                                                                       \right] \in \R^{Nn (n+m+p+q)}.
\end{equation*}
Also, let $\bar{A}_i=A+\delta_i A$, $\bar{B}_i=B+\delta_i B,$ $\bar{E}_i=E+\delta_i E$, and $\bar{C}_i=C+\delta_i C$. Then \eqref{system0} can be put in the following form:
\begin{equation}\label{system}
\begin{split}
x_i(t+1)=&\bar{A}_ix_i(t)+\bar{B}_iu_i(t-r_{con})+\bar{E}_iv(t),\\
y_i(t)=&\bar{C}_ix_i(t).\\
\end{split}
\end{equation}
The plant \eqref{system} and \eqref{exo} can be viewed as a multi-agent system with the exosystem \eqref{exo} as the leader and the $N$ subsystems of \eqref{system} as the followers. The communication
topology can be described by a directed graph $\mathcal{\bar{G}}=(\mathcal{\bar{V}},\mathcal{\bar{E}})$, where $\mathcal{\bar{V}}=\{0,1,\dots,N\}$ is the node set with the node 0 associated with the exosystem \eqref{exo} and all the other nodes associated with the $N$ subsystems \eqref{system}, and $\mathcal{\bar{E}}$ is the edge set. The edge $(j,i) \in \mathcal{\bar{E}},\ i\neq j,\ i,j=0,\dots,N$, if and only if the control $u_i,\ i=1,\dots,N,$ can access the state $x_j$ and/or the output $y_j$ of subsystem $j,\ j=0,\dots,N$.
If $(j,i) \in \mathcal{\bar{E}}$, node $j$ is called a neighbor of the node $i$. We use $\bar{\mathcal{{N}}}_i$ to denote the neighbor set of node $i$ with respect to $\mathcal{\bar{V}}$.

Due to the communication constraint, we are limited to consider the class of distributed control laws. Mathematically, such a control law is described as follows,

\begin{equation}\label{con}
\begin{split}
&u_i(t)\!=\!k_i(z_i(t), z_j(t), x_i(t\!-\!r_{com})\!-\!x_j(t\!-\!r_{com}), j\in \bar{\mathcal{N}}_i),\\
&z_i(t+1)\!=\!g_i(z_i(t), y_i(t\!-\!r_{com})\!-\!y_j(t\!-\!r_{com}), j\in \bar{\mathcal{N}}_i),
\end{split}
\end{equation}
where $z_i\in\R^{n_{zi}}$, $k_i$ and $g_i$ are linear functions of their arguments, $r_{com}\in \mathbb{Z}^{+}$ represents the communication delay among agents. The control law \eqref{con} is called a distributed dynamic state feedback control law, and is further called dynamic output feedback control law if the function $k_i$ is independent of any state variable.

Now, we can state our problem as follows:
\begin{Definition}\label{Def1.1}
Discrete-time linear cooperative robust output regulation problem: given the multi-agent system \eqref{system}, the exosystem \eqref{exo}, and a digraph $\mathcal{\bar{G}}$, design a control law of the form \eqref{con} such that the closed-loop system satisfies the following properties.
\begin{Property}\label{Per1}
The nominal closed-loop system is exponentially stable when $v=0$.
\end{Property}

\begin{Property}\label{Per2}
There exists an open neighborhood $W$ of $w=0$ such that, for any $w \in W$ and any initial conditions $x_{i0}$, $z_{i0}$ and $v_0$, the regulated output $\lim_{t \rightarrow \infty}e_i(t)= 0,\ i=1,\dots,N$.
\end{Property}
\end{Definition}

\vspace{-6pt}

\section{A GENERAL FRAMEWORK}\label{frame1}
\vspace{-2pt}

It is known that, the robust output regulation problem of a delay-free plant can be converted to the stabilization problem of an augmented system composed of the given plant and a dynamic compensator called internal model \cite{Da, Fr, FW}. This design philosophy is known as the  internal model principle. Paper \cite{JSSC} has generalized the internal model design from delay-free discrete-time systems to discrete-time systems with both input and communication delays. In this section, we will further generalize the framework in \cite{JSSC} for a single system to multi-agent systems. This framework will be based on the concept of the distributed internal model. For this purpose, we will first recall the concept of the minimal $p$-copy internal model as follows.
\begin{Definition}\label{Definter}
A pair of matrices $(G_1,G_2)$ is said to be a minimal $p$-copy
internal model of the matrix $S$ if the pair takes the following
form:
\begin{equation}\label{inter}
G_1= I_p \otimes \beta,\ \
G_2= I_p \otimes \sigma,\\
\end{equation}
where $\beta$ is a constant square matrix whose characteristic polynomial equals the minimal polynomial of $S$, and $\sigma$ is a constant column vector such that $(\beta,\sigma)$ is controllable.
\end{Definition}


To introduce the distributed internal model, let $\mathcal{\bar{A}}=\left[ a_{ij}\right] \in R^{(N+1)\times(N+1)}$ and $\mathcal{\bar{L}}=[l_{ij}] \in R^{(N+1)\times(N+1)}$ be the weighted adjacency matrix and Laplacian of the digraph $\mathcal{\bar{G}}$,  respectively.
In terms of the elements of $\mathcal{\bar{A}}$, we can define a virtual regulated output $e_{vi}(t)$ for each follower subsystem $i$ as follows:
\begin{equation}\label{sysvirt}
e_{vi}(t) =\sum_{j\in \bar{\mathcal{{N}}}_i} a_{ij} (y_i(t)-y_j(t)),\ i=1,\dots,N.
\end{equation}
It is noted that the subsystem $e_{vi}(t)$ can access the regulated error $(y_i(t)-y_j(t))$ if and only if the node $j$ is the neighbor of the node $i$.

We call the following dynamic compensator
\begin{equation}\label{internal}
z_i(t+1)=G_1z_i(t)+G_2e_{vi}(t-r_{com}), i = 1, \cdots, N
\end{equation}
a distributed internal model of the plant \eqref{system} and the exosystem \eqref{exo}.

\begin{Remark}\label{Remeig}
Let $e=\mbox{col}(e_1,\dots,e_{N})$ and $e_v=\mbox{col}(e_{v1},\dots,e_{vN})$. Then it can be verified that $e_v = (H \otimes I_p) e$, where $H \in \R^{N\times N}$ consists of the last $N$ rows and the last $N$ columns of $\mathcal{\bar{L}}$.  By Lemma 1 of \cite{Su1}, the matrix $-H$ is Hurwitz if and only if the digraph is connected. Thus, if the digraph is connected, then $e = 0$ iff $e_v = 0$.
\end{Remark}


Having introduced the $p$-copy internal model and defined the virtual regulated output $e_{vi}(t)$, we can describe our control laws as follows:

\noindent 1) Distributed dynamic state feedback control law
\begin{equation}\label{ctr1}
\begin{split}
u_i(t)=&  K_x \eta_i(t) + K_z z_i(t),\\
z_i(t+1) =&  {G}_1 z_i(t) + {G}_2 e_{vi}(t-r_{com}),\ i=1,\dots,N,  \\
\end{split}
\end{equation}
where $\eta_i(t)= \sum_{j\in \bar{\mathcal{{N}}}_i} a_{ij}\left(x_i\left(t-r_{com}\right)-x_j\left(t-r_{com}\right)\right) $, $x_{0}(t)=0$, $z_i(t) \in \R^{n_{zi}}$ with $n_{zi}$ to be specified later, $(K_x,K_z)$ are constant matrices of appropriate dimensions to be designed later, $(G_1,G_2)$ are defined in \eqref{inter}.

\noindent 2) Distributed dynamic output feedback control law
\begin{equation}\label{ctr2}
\begin{aligned}
u_i(t)=&K_1 z_i(t) + K_{2} \hat{\eta}_i(t), \\
z_i(t+1)\! =&  {G}_1 z_i(t) + {G}_2 e_{vi}(t-r_{com}),  \\
\xi_i(t+1)\! =& A \xi_i(t)\! +\! Bu_i(t-r) \! -\!LC \hat{\eta}_i(t)\! +\!Le_{vi}(t-r_{com}), \ i=1,\dots,N,
\end{aligned}
\end{equation}
where $\xi_i(t) \in \R^{n_i}$, $\hat{\eta}_i(t)= \sum_{j\in \bar{\mathcal{{N}}}_i}a_{ij}\left(\xi_i(t)-\xi_j(t)\right) $, $\xi_0(t)=0$, and $z_i(t) \in \R^{n_{zi}}$ with $n_{zi}$ to be specified later, $r=r_{com}+r_{con}$, $(K_1,K_{2},L)$ are constant matrices of appropriate dimensions to be designed later and $(G_1,G_2)$ are defined in \eqref{inter}.

\begin{Remark}\label{Remtrans}
To handle the communication delay, introduce the coordinate transformation $z_i(t)=\bar{z}_i(t-r_{com}),$ and $\xi_i(t)=\bar{\xi}_i(t-r_{com})$. Then,  the state feedback control law \eqref{ctr1} becomes as follows:
\begin{equation}\label{ctr1tr}
\begin{split}
u_i(t)=&  K_x \eta_i(t) + K_z \bar{z}_i(t-r_{com}),\\
\bar{z}_i(t+1) =&  {G}_1 \bar{z}_i(t) + {G}_2 e_{vi}(t) ,\ i=1,\dots,N,
\end{split}
\end{equation}
and, respectively, the output feedback control law \eqref{ctr2} becomes as follows:
\begin{equation}\label{ctr22}
\begin{split}
u_i(t)=& K_1 \bar{z}_i(t-r_{com})  +  K_2 \bar{\eta}_i(t-r_{com}),\\
\bar{z}_i(t+1) =&  {G}_1 \bar{z}_i(t) + {G}_2 e_{vi}(t),\\
\bar{\xi}_i(t+1) =& A \bar{\xi}_i(t) + Bu_i(t-r_{con}) -LC \bar{\eta}_i(t) + Le_{vi}(t),
\ i=1,\dots,N,
\end{split}
\end{equation}
where $\bar{\eta}_i(t)=\sum_{j\in\bar{\mathcal{{N}}}_i}a_{ij}\left(\bar{\xi}_i(t)-\bar{\xi}_j(t)\right)$.
\end{Remark}

Attaching the distributed internal model \eqref{internal} to the state equation of the plant \eqref{system} leads to the following so-called the auxiliary augmented system of \eqref{system}:
\begin{equation}\label{syscomp}
\begin{split}
x_i(t+1) &= \bar{A}_i x_i(t) +  \bar{B}_i u_i(t-r_{con}) + \bar{E}_i v(t), \\
z_i(t+1) & =  {G}_1 z_i (t) +  {G}_2 e_{vi} (t-r_{com}),~ i = 1, \cdots, N. 
\end{split}
\end{equation}


We now ready to present our main result of this section as follows:

\begin{Lemma}\label{lem3.1}
Suppose $S$ has no eigenvalues with modulus smaller than $1$ and the digraph $\mathcal{\bar{G}}$  is connected. Then,

(i) if a static state feedback control law of the form
\begin{equation}\label{Ctr1n10}
u_i(t)= K_x \eta_i(t) + K_z \bar{z}_i(t-r_{com}),\ i=1,\dots,N
\end{equation}
stabilizes the nominal plant of the auxiliary  augmented system \eqref{syscomp} with $v=0$, then, the  dynamic state feedback control law
\eqref{ctr1tr}
solves the cooperative robust output regulation problem of the  plant \eqref{system} and the exosystem \eqref{exo}.

(ii) if a dynamic output feedback control law of the form
\begin{equation} \label{Ctr1n20}
\begin{split}
u_i(t)=& K_1 \bar{z}_i(t-r_{com})  +  K_2 \bar{\eta}_i(t-r_{com}),\\
\bar{\xi}_i(t+1) =& A \bar{\xi}_i(t) + Bu_i(t-r_{con}) -LC \bar{\eta}_i(t) + Le_{vi}(t),
\ i=1,\dots,N,
\end{split}
\end{equation}
 stabilizes the nominal plant of the auxiliary  augmented system \eqref{syscomp} with $v=0$, then, the  dynamic output feedback control law
\eqref{ctr22}
solves the cooperative robust output regulation problem of the  plant \eqref{system} and the exosystem \eqref{exo}.
\end{Lemma}

\begin{Proof}
Let $x=\mbox{col}(x_1,\dots,x_N)$,  $z= \mbox{col}(z_1,\dots,z_N)$, $\xi= \mbox{col}(\xi_1,\dots,\xi_N)$, $u=\mbox{col}(u_1,\dots,u_N)$, $\bar{A}=\text{block~diag}(\bar{A}_1,$ $\dots,\bar{A}_N)$, $\bar{B}=\text{block~diag}(\bar{B}_1,\dots,\bar{B}_N)$, $\bar{E}=(\bar{E}_1^T,\dots,$ $\bar{E}_N^T)^T$,  $\bar{G}_1=I_N\otimes G_1$ and $\bar{G}_2=I_N\otimes G_2$. Then 
the auxiliary augmented system  \eqref{syscomp} can be put into the following compact form:
\begin{equation}\label{Csyscomp}
\begin{split}
x(t+1) &= \bar{A} x(t) +  \bar{B} u(t-r_{con}) + \bar{E} v(t), \\
z(t+1) & =  \bar{ {G}}_1 z (t) +  \bar{ {G}}_2 e_v(t-r_{com}).\\
\end{split}
\end{equation}

Let $\bar{C}=(H\otimes I_p)~\text{block~diag}~(\bar{C}_1,\dots,\bar{C}_N)$, $\bar{F}=(\Delta \textbf{1}_N)\otimes F$. Then the virtual regulated output can be put in the following compact form:
\begin{equation}\label{vertualoutputev}
e_v(t) = \bar{C} x(t) + \bar{F} v(t).
\end{equation} 
 
Define a so-called auxiliary system as follows:
\begin{equation}\label{Csystemauxi}
\begin{split}
x(t+1) &= \bar{A} x(t) +  \bar{B} u(t-r_{con}) + \bar{E} v(t),\\
v(t+1) &= S v(t),\\
e_v(t) &= \bar{C} x(t) + \bar{F} v(t).\\
\end{split}
\end{equation}

Further, let $\bar{K}_x=H\otimes K_x$, $\bar{K}_z=I_N\otimes K_z$,
$\bar{K}_1=I_N \otimes K_1$, $\bar{K}_2=H\otimes K_2$, $\bar{S}_1=I_N\otimes A-H\otimes LC$, $\bar{S}_2=I_N\otimes L$, $\bar{S}_3=\bar{B}\left(\begin{array}{cc}\bar{K}_1&\bar{K}_2\\ \end{array}\right)$, and $\bar{\zeta}=\mbox{col}(\bar{z},\bar{\xi})$. Then,  the static state feedback control law  (\ref{Ctr1n10}) and the dynamic output feedback control law (\ref{Ctr1n20}) can be put into the following compact form:
\begin{equation}\label{sCtr1n1}
\begin{split}
u(t)=& \bar{K}_x x(t-r_{com})+ \bar{K}_z \bar{z} (t-r_{com}),
\end{split}
\end{equation}
and, respectively,
\begin{equation}\label{sCtr1n2}
\begin{split}
u(t)=& \bar{K}_1 \bar{z}(t-r_{com}) + \bar{K}_{2} \bar{\xi}(t-r_{com}),\\
\bar{\xi}(t+1)=& \bar{S}_{1}\bar{\xi}(t)+ \bar{S}_2 e_v(t) + \bar{S}_{3} \bar{\zeta}(t- r).
\end{split}
\end{equation}
Similarly, the dynamic state feedback control law  \eqref{ctr1tr}  and the dynamic output feedback control law \eqref{ctr22} can be put into the following compact form:
\begin{equation}\label{Ctr1n1}
\begin{split}
u(t)=& \bar{K}_x x(t-r_{com})+ \bar{K}_z \bar{z} (t-r_{com}),\\
\bar{z}(t+1) =& {\bar G}_1 \bar{z} (t) +{\bar G}_2 e_v(t),  \\
\end{split}
\end{equation}
and, respectively,
\begin{equation}\label{Ctr1n2}
\begin{split}
u(t)=& \bar{K}_1 \bar{z}(t-r_{com}) + \bar{K}_{2} \bar{\xi}(t-r_{com}),\\
\bar{z}(t+1) =& \bar{ {G}}_1 \bar{z} (t) +  \bar{ {G}}_2 e_v(t), \\
\bar{\xi}(t+1)=& \bar{S}_{1}\bar{\xi}(t)+ \bar{S}_2 e_v(t) + \bar{S}_{3} \bar{\zeta}(t- r).
\end{split}
\end{equation}

Now applying Lemma 3.1 of  \cite{JSSC} to the auxiliary system (\ref{Csystemauxi})  viewing $e_v$ as the tracking error concludes that
if the static state feedback control law (\ref{sCtr1n1}) or, respectively,  the dynamic output feedback control law (\ref{sCtr1n2})
stabilizes the nominal plant of the auxiliary  augmented system \eqref{Csyscomp} with $v=0$, then, the  dynamic state feedback control law
(\ref{Ctr1n1}), or, respectively, the dynamic output feedback control law (\ref{Ctr1n2})
solves the robust output regulation problem of the the auxiliary system (\ref{Csystemauxi}).

Finally, by Remark \ref{Remeig}, if the digraph $\mathcal{\bar{G}}$ is connected, then $e = 0$ iff $e_v = 0$. Thus, under the assumption that the digraph $\mathcal{\bar{G}}$ is connected, the control law \eqref{Ctr1n1} or the control law \eqref{Ctr1n2} solves the robust output regulation problem of the auxiliary system (\ref{Csystemauxi})  viewing $e_v$ as the tracking error
 iff the same control law
solves the cooperative robust output regulation problem of the plant \eqref{system} and the exosystem \eqref{exo}.
\end{Proof}

\section{MAIN RESULT}\label{result}

In this section, we will present the main results of the cooperative robust output regulation problem based on the internal model framework introduced in section \ref{frame1}. By Lemma \ref{lem3.1}, it suffices to
 stabilize the auxiliary augmented system \eqref{syscomp} by either distributed dynamic state feedback control law \eqref{Ctr1n10}
 or distributed dynamic output feedback control law \eqref{Ctr1n20}. Before presenting our main result, we need the following assumptions.


\begin{Assumption}\label{Ass2.3}
The matrix pair $(A,B)$ is stabilizable.
\end{Assumption}

\begin{Assumption}\label{Ass2.4}
The matrix pair $(C,A)$ is detectable.
\end{Assumption}

\begin{Assumption}\label{Ass2.5}
For all $\lambda \in \sigma(S)$, 
\begin{equation}\label{tzero}
\text{rank} \left(
                                      \begin{array}{cc}
                                          A -\lambda I_n &  B \\
                                          C & 0 \\
                                      \end{array}
                                    \right)=n+p.
\end{equation}
\end{Assumption}

\begin{Assumption}\label{Ass2.6}
The digraph $\bar{\mathcal{G}}$ contains a directed spanning tree with the node $0$ as the root.
\end{Assumption}

\begin{Assumption}\label{Ass2.2}
All the eigenvalues of $S$ are on the unit circle.
\end{Assumption}

\begin{Assumption}\label{Ass2.7}
$A$ has no eigenvalues with modulus greater than $1$.
\end{Assumption}
\begin{Remark}\label{afterassu}
Assumptions \ref{Ass2.3} to \ref{Ass2.6} are quite standard  and they are also needed in \cite{Su3} for the cooperative output regulation problem of continuous-time systems even if there are no communication delay and input delay. Assumptions \ref{Ass2.2} and \ref{Ass2.7} are additional and they are made so that the delayed system can be stabilized by using the method in \cite{Bzhoubook}, which is summarized in Lemma \ref{Lem4.1ap} and Lemma \ref{Lem4.3ap} of the Appendix. These two assumptions can be removed if there are no communication and input delays.
\end{Remark}

Now we establish some lemmas to lay the foundation for our main results.

\begin{Lemma} \label{Lem4.3x}
Let  ${A} \in \R^{n \times n}$, ${B} \in \R^{n \times m}$, and $H \in \R^{N \times N}$. Suppose all the  eigenvalues of $A$ have modulus equal to or smaller than $1$,
all the eigenvalues of $H$ have positive real parts, and $(A, B)$ is stablizable. Then, there exists a matrix $K \in \R^{m \times n}$ such that the matrix $\left(I_N \otimes A + H \otimes (BK)\right)$ is Schur.
\end{Lemma}

\begin{Proof}
Denote the eigenvalues of $H$ by $\lambda_i,\ i=1,\dots,N$ where, for $i=1,\dots,N$,   $\lambda_i$ has positive real part by assumption.
Let $T_1$ be the non-singular matrix such that $J_H=T_1 H T_1^{-1}$ is a lower triangular matrix with its $i^{th}$ diagonal elements being denoted by $\lambda_i$. Then $(T_1 \otimes I_N) (I_N \otimes A + H \otimes (BK))(T_1 \otimes I_N)^{-1} =
(I_N \otimes A + J_H \otimes (BK))$ is a lower triangular system whose diagonal blocks are of the form $A+\lambda_i BK$, $i = 1, \cdots, N$.
Now define the following systems
\begin{equation}\label{sta1}
x_i (t+1) = A x_i (t) + \lambda_i B u_i(t), ~ i = 1, \cdots, N.
\end{equation}
Then,  by Lemma \ref{Lem4.3ap} of the Appendix, there exists a  matrix $K \in \R^{m \times n}$ such that $A+\lambda_i BK$, $i = 1, \cdots, N$, are Schur. The proof is completed.
\end{Proof}

\begin{Lemma}\label{Lem4.4}
Consider the system of the form
\begin{equation}\label{linm1}
\begin{aligned}
x_c(t+1) =& \left(
               \begin{array}{cc}
                  I_N\otimes A  & 0 \\
                  H\otimes G_2C  &  I_N\otimes G_1  \\
               \end{array}
             \right) x_c(t)  +\left(
                              \begin{array}{cc}
                                H\otimes B & I_N\otimes B \\
                                0_{N n_z \times Nm} & 0_{N n_z \times Nm} \\
                              \end{array}
                            \right)u_c(t-r),\\
\end{aligned}
\end{equation}
where $x_c \in \R^{N(n+n_z)}$, $u_c \in \R^{2Nm}$,  $(G_1,G_2)$ is the minimal p-copy internal model of $S$ as defined in \eqref{inter},
and $x_{c0} \in \mathcal{C}\big(I[-r,0]$ $, \R^{N(n+n_z)}\big)$. Then, under Assumptions \ref{Ass2.3}, \ref{Ass2.5}, \ref{Ass2.2}, and \ref{Ass2.7}, there exist matrices $K_x \in\R^{m \times n}$ and $K_z \in \R^{m \times n_z}$, such that under the state feedback control law
$u_c(t)=\left(
          \begin{array}{cc}
            I_N \otimes K_x & 0_{Nm \times N n_z} \\
            0_{Nm \times N n} & I_N \otimes K_z \\
          \end{array}
        \right)x_c(t)$, system \eqref{linm1} is asymptotically stable if and only if Assumption \ref{Ass2.6} is satisfied.
\end{Lemma}
\begin{Proof}
 This lemma can be viewed as a discrete-time counterpart of Lemma 4.3 of \cite{Lu2015}, and its proof is also similar to that of Lemma 4.3 of \cite{Lu2015}.
 In particular, the proof of the necessary part is the same as the proof of the necessary part of Lemma 4.3 of \cite{Lu2015}.
 Thus we only focus on the sufficient part.

 As in Lemma \ref{Lem4.3x}, denote the eigenvalues of $H$ by $\lambda_i,\ i=1,\dots,N$. From the proof of Lemma 4.3 of \cite{Lu2015}, there exist nonsingular matrices $T_x \in \R^{N(n+n_z) \times N(n+n_z)}$ and $T_u \in \R^{2Nm \times 2Nm}$ such that $\tilde{x}_c = T_x x_c$ with input $\tilde{u}_c = T_u u_c$ is governed by
\begin{equation}\label{linm4}
\begin{aligned}
\tilde{x}_{ci}(t+1) =& \left(
               \begin{array}{cc}
                  A & 0 \\
                  G_2C & G_1 \\
               \end{array}
             \right) \tilde{x}_{ci}(t) + \lambda_i
             \left(
               \begin{array}{cc}
                 B & B \\
                 0 & 0 \\
               \end{array}
             \right)\tilde{u}_{ci} (t-r),~ i = 1, \cdots, N, \\
\end{aligned}
\end{equation}
where $\tilde{x}_{c} =\mbox{col}(\tilde{x}_{c1},\dots,\tilde{x}_{cN})$ with $\tilde{x}_{ci} \in \R^{(n+n_z)}$ and
$\tilde{u}_{c} =\mbox{col}(\tilde{u}_{c1},\dots,\tilde{u}_{cN})$ with $\tilde{u}_{ci} \in \R^m$, and
\begin{equation} \label{linma}
{T}_u^{-1} \left(
                     \begin{array}{cc}
                       I_N\otimes K_x & 0 \\
                       0 & I_N\otimes K_z \\
                     \end{array}
                   \right) T_x = \left(
                     \begin{array}{cc}
                       I_N\otimes K_x & 0 \\
                       0 & I_N\otimes K_z \\
                       \end{array}
                   \right).
\end{equation}
Let $\mathcal{A} = \left(
               \begin{array}{cc}
                  A & 0 \\
                  G_2C & G_1 \\
               \end{array}
             \right)$ and $\mathcal{B} = \left(
                               \begin{array}{c}
                                 B \\
                                 0 \\
                               \end{array}
                             \right)$.
Under Assumptions \ref{Ass2.3}, \ref{Ass2.5} and \ref{Ass2.2}, by Lemma 1.37 of \cite{J.H}, $(\mathcal{A},\mathcal{B})$ is stabilizable. Moreover, under additional Assumption \ref{Ass2.7},  $\mathcal{A}$ has no eigenvalues with modulus greater than $1$.
By Lemma \ref{Lem4.3ap}, there exists a matrix $\tilde{K}=(K_x,K_z)$, where $K_x \in \R^{m \times n}$ and $K_z \in \R^{m \times n_z}$ such that
the following systems
\begin{equation}\label{linm411}
\begin{aligned}
\tilde{x}_{ci}(t+1) =& \mathcal{A}   \tilde{x}_{ci}(t)
                            +
                             \lambda_i \mathcal{B} \tilde{K} \tilde{x}_{ci}(t-r), \ i=1,\dots,N,
\end{aligned}
\end{equation}
are asymptotically stable.  Thus, for each $i = 1, \cdots, N$, the state feedback control law
 $ \tilde{u}_{ci} (t) =  \left(
                                          \begin{array}{cc}
                                            K_x & 0 \\
                                            0 & K_z \\
                                          \end{array}
                                        \right)
 \tilde{x}_{ci}(t) $  asymptotically stabilizes the system  (\ref{linm4}).

 Finally, from (\ref{linma}), we have
\begin{equation}
\begin{aligned}
{u}_{c}  &=  T^{-1}_u \tilde{u}_c \\
& =  T^{-1}_u \left(
                     \begin{array}{cc}
                       I_N\otimes K_x & 0 \\
                       0 & I_N\otimes K_z \\
                     \end{array}
                   \right)\tilde{x}_c \\
                   &=  T^{-1}_u \left(
                     \begin{array}{cc}
                       I_N\otimes K_x & 0 \\
                       0 & I_N\otimes K_z \\
                     \end{array}
                   \right) T_x  {x}_c \\
                   &
                   = \left(
                     \begin{array}{cc}
                       I_N\otimes K_x & 0 \\
                       0 & I_N\otimes K_z \\
                       \end{array}
                   \right) {x}_c.
\end{aligned}
\end{equation}
Thus, the proof is completed.

\end{Proof}

\begin{Theorem}\label{Theo4.1}
Under Assumptions \ref{Ass2.3}, \ref{Ass2.5}, \ref{Ass2.2} and \ref{Ass2.7}, there exist matrices $K_x \in \R^{m \times n}$, $K_z \in \R^{m \times n_z}$ such that the cooperative robust output regulation problem is solved by the distributed dynamic state
feedback control law \eqref{ctr1tr} with $(G_1,G_2)$ being the minimal $p$-copy internal model of $S$ if and only if Assumption \ref{Ass2.6} is satisfied.
\end{Theorem}
\begin{Proof}
Let $x_c=\mbox{col}(x,\bar{z})$. Then the closed-loop system composed of the auxiliary  system \eqref{Csystemauxi} and dynamic state feedback control law \eqref{Ctr1n1} is the same as
the closed-loop system composed of the auxiliary augmented system \eqref{syscomp} and the static state feedback control law \eqref{Ctr1n10} and
can be put into the following form:
\begin{equation}\label{cl3}
\begin{split}
x_c(t+1) &= \sum_{l=0}^{1} A_{cwl} x_c(t-\bar{r}_l) + B_{cw} v(t),\\
e_v(t) &=  C_{cw} x_c(t) + D_{cw} v(t), \\
\end{split}
\end{equation}
where $\bar{r}_0=0$, $\bar{r}_1=r$, and
\begin{equation*}
\begin{split}
A_{cw0} =& \left(
          \begin{array}{cc}
            \bar{A}  & 0 \\
            \bar{G}_2 \bar{C}  & \bar{G}_1 \\
          \end{array}
        \right), \
A_{cw1} =  \left(
          \begin{array}{cc}
            \bar{B} \bar{K}_x & \bar{B} \bar{K}_z \\
            0 & 0  \\
          \end{array}
        \right), \\
B_{cw}= & \left(
            \begin{array}{c}
              \bar{E} \\
            \bar{G}_2 \bar{F}  \\
            \end{array}
          \right) , \
C_{cw} =  \left(
             \begin{array}{cc}
               \bar{C} & 0\\
             \end{array}
           \right),
D_{cw}= \bar{F}. \\
\end{split}
\end{equation*}
Thus, the  nominal closed-loop system with $v$ set to $0$ is as follows:
\begin{equation}\label{eqthe1}
\begin{aligned}
x_c(t+1) =& \left(
               \begin{array}{cc}
                  I_N\otimes A  & 0 \\
                  H\otimes G_2C  &  I_N\otimes G_1  \\
               \end{array}
             \right) x_c(t) +
             \left(
               \begin{array}{cc}
                H\otimes B & I_N\otimes B \\
                 0_{Nn_z \times Nm} & 0_{Nn_z \times Nm} \\
               \end{array}
             \right)\\ \times
             &\left(
               \begin{array}{cc}
                 I_N \otimes K_x & 0_{Nm \times Nn_z} \\
                 0_{Nm \times Nn} & I_N \otimes K_z \\
               \end{array}
             \right) x_c(t-r).
\end{aligned}
\end{equation}
By Lemma \ref{Lem4.4}, there exist matrices $K_x \in \R^{m \times n}$ and $K_z \in \R^{m \times n_z}$, such that system \eqref{eqthe1} is asymptotically stable. The proof is thus completed by
invoking Lemma \ref{lem3.1}.  \end{Proof}

To study the output feedback case, we need the following lemma.
\begin{Lemma}\label{Lem4.5}
Consider the system of the form
\begin{equation}\label{Thee1}
\begin{aligned}
x_c(t+1) \!=& \!\left(\!\!
               \begin{array}{ccc}
                  I_N\otimes A  \!& \!0 \!&\! 0 \\
                  H\otimes G_2 C \! &\! I_N \otimes G_1\! &\! 0 \\
                  H\otimes LC  \!&\!  0  \!&\! I_N\otimes A - H\otimes LC\\
               \end{array}\!\!
             \right)x_c(t) \! \\&+ \!\left(\!\!
                              \begin{array}{ccc}
                                0\!& \! I_N\otimes B \! & \! H\otimes B  \\
                                0\!&\! 0\! &\! 0\\
                                0\! &\! I_N\otimes B \!& \!H\otimes B \\
                              \end{array}\!\!
                            \right)\!\!u_c(t-r),\\
\end{aligned}
\end{equation}
where $x_c \in \R^{N(2n+n_z)}$, $u_c \in \R^{3Nm}$, $(G_1,G_2)$ is the minimal p-copy internal model of $S$ as defined in \eqref{inter}, and
$x_{c0} \in \mathcal{C}\big(I[-r,0]$ $, \R^{N(2n+n_z)}\big)$. Then, under Assumptions \ref{Ass2.3}-\ref{Ass2.5}, \ref{Ass2.2} and \ref{Ass2.7}, there exist matrices $K_1 \in \R^{m \times n_z}$, $K_2 \in \R^{m \times n}$ and $L \in \R^{n \times p}$, such that under the state feedback control law $u_c(t)=Kx_c(t)$, where $K=\left(
                                                                                  \begin{array}{ccc}
                                                                                    I_N\otimes K_2 & 0 & 0 \\
                                                                                    0 & I_N\otimes K_1 & 0 \\
                                                                                    0 & 0 & I_N\otimes K_2 \\
                                                                                  \end{array}
                                                                                \right)$, system \eqref{Thee1} is asymptotically stable if and only if Assumption \ref{Ass2.6} is satisfied.
\end{Lemma}
\begin{Proof}
This lemma can be viewed as a discrete-time counterpart of Lemma 4.4 of \cite{Lu2015}, and its proof is also similar to that of Lemma 4.4 of \cite{Lu2015}.

Let
$T_x=\left(
     \begin{array}{ccc}
       I_{Nn} & 0 & 0 \\
       0 & I_{Nn_z} & 0 \\
       -I_{Nn} & 0 & I_{Nn} \\
     \end{array}
   \right)
$ and ${T}_u=\left(
                                         \begin{array}{ccc}
                                           I_{Nm} & 0 & 0 \\
                                           0 & I_{Nm}& 0  \\
                                           -I_{Nm} & 0 & I_{Nm} \\
                                         \end{array}
                                       \right)
$. Then the state $\bar{x}_c = T_x x_c$ with input $\bar{u}_c = T_u u_c$ is governed by
\begin{equation}\label{tThee2}
\begin{aligned}
\bar{x}_c(t+1)=& \left(
               \begin{array}{ccc}
                  I_N\otimes A \!\! & 0\!\! & 0 \\
                  H\otimes G_2 C \!\! & I_N \otimes G_1 \!\!& 0 \\
                  0\!\! &  0  & \!\! I_N\otimes A - H\otimes LC\\
               \end{array}
           \right)\bar{x}_c(t)   \\&+  \left(\!
                              \begin{array}{ccc}
                                 H\otimes B \!&  I_N\otimes B  \!& H\otimes B   \\
                                0 \!& 0\!& 0\\
                                0\!& 0 \!& 0\\
                              \end{array}
                          \! \right)   \bar{u}_c(t-r).
\end{aligned}
\end{equation}

Denote $\bar{x}_c=\mbox{col}(\bar{x}_{c1},\bar{x}_{c2})$ with $\bar{x}_{c1} \in \R^{N(n+n_z)}$ and $\bar{x}_{c2} \in \R^{Nn}$.
Then, by Lemma \ref{Lem4.4}, under Assumptions \ref{Ass2.3}, \ref{Ass2.5}, \ref{Ass2.2} and \ref{Ass2.7}, there exist matrices $K_1 \in \R^{m \times n_z}$ and $K_2 \in \R^{m \times n}$, such that
 the following system
\begin{equation}\label{Thee2}
\begin{aligned}
\bar{x}_{c1}(t+1) =& \left(
               \begin{array}{cc}
                  I_N\otimes A  & 0  \\
                  H\otimes G_2 C  & I_N \otimes G_1 \\
               \end{array}
             \right) \bar{x}_{c1}(t)   + \left(
                              \begin{array}{cc}
                                 H\otimes B&  I_N\otimes B\\
                                0 & 0\\
                              \end{array}
                            \right) \\&\times   \left(
                                       \begin{array}{cc}
                                         I_N \otimes  K_2& 0 \\
                                         0 &  I_N \otimes  K_1 \\
                                       \end{array}
                                     \right)
                            \bar{x}_{c1}(t-r)\\
\end{aligned}
\end{equation}
is asymptotically stable if and only if the digraph satisfies Assumption \ref{Ass2.6}.

Let $\bar{K}=\left(
                                                \begin{array}{ccc}
                                                  I_N\otimes K_2 & 0 & 0 \\
                                                  0 & I_N\otimes K_1 & 0 \\
                                                  0 & 0 & I_N\otimes K_2 \\
                                                \end{array} \right)$.
Then, under the state feedback control law $ u_c(t)= {T}^{-1}_u \bar{K}\bar{x}_c(t)$, the closed-loop system of \eqref{tThee2} is as follows:
\begin{equation}\label{Eqlemc1}
\begin{split}
\bar{x}_{c1}(t+1) =& \left(
               \begin{array}{cc}
                  I_N\otimes A  & 0  \\
                  H\otimes G_2 C  & I_N \otimes G_1 \\
               \end{array}
             \right) \bar{x}_{c1}(t)   + \left(
                              \begin{array}{cc}
                                 H\otimes B&  I_N\otimes B\\
                                0 & 0\\
                              \end{array}
                            \right) \\&\times   \left(
                                       \begin{array}{cc}
                                         I_N \otimes  K_2& 0 \\
                                         0 &  I_N \otimes  K_1 \\
                                       \end{array}
                                     \right)
                            \bar{x}_{c1}(t-r)  + \left(
                                                         \begin{array}{c}
                                                           H\otimes B \\
                                                           0 \\
                                                         \end{array}
                                                       \right) (I_N \otimes  K_2) \bar{x}_{c2}(t-r),\\
\bar{x}_{c2}(t+1) = &(I_N\otimes A - H\otimes LC) \bar{x}_{c2}(t).
\end{split}
\end{equation}

Note that, $(I_N\otimes A - H\otimes LC)^T = (I_N\otimes A^T - H^T \otimes C^T L^T)$. Under the assumptions of this lemma, $(A^T, C^T)$ is stablizable, all the eigenvalues of $A^T$ have modulus equal to or smaller than $1$, and all the eigenvalues of $H^T$ have positive real parts. By Lemma \ref{Lem4.3x}, there exists a matrix $K \in \R^{m \times n}$ such that the matrix $(I_N\otimes A^T + H^T \otimes C^T K)$ is Schur, which implies,  with $L = - K^T$, $\left(I_N\otimes A - H\otimes LC\right)$ is Schur. Moreover by Lemma \ref{Lem4.4}, the $ \bar{x}_{c1}$ subsystem with $\bar{x}_{c2}$ set to zero is asymptotically stable. Thus, by Lemma 3 of \cite{IET}, system \eqref{Eqlemc1} is asymptotically stable.
Furthermore, since $\bar{x}_c(t)=Tx_c(t)$, we have
\begin{equation}\label{}
\begin{split}
&u_c(t)\\=&{T}^{-1}_u \bar{K}\bar{x}_c(t) \\=& {T}^{-1}_u \bar{K} T_x x_c(t)\\
=& \left(\!\!\!
     \begin{array}{ccc}
       I_{Nm} \!& \!0 \!& \!0 \\
       0 \!& \!I_{Nm} \!& \!0 \\
       I_{Nm} \!& \!0\! &\! I_{Nm} \\
     \end{array}\!\!\!
   \right)\!\! \left(\!\!\!
            \begin{array}{ccc}
              I_N\otimes K_2 \!&\! 0 \!&\! 0 \\
              0 \!&\! I_N\otimes K_1 \!&\! 0 \\
               0 \!&\! 0 \!&\! I_N\otimes K_2 \\
            \end{array}\!\!\!
          \right)
          \left(
            \begin{array}{ccc}
              I_{Nn} & 0 & 0 \\
              0 & I_{Nn_z} & 0 \\
              -I_{Nn} &  0 & I_{Nn}\\
            \end{array}
          \right) x_c(t)\\
=&\left(
     \begin{array}{ccc}
       I_N\otimes K_2 & 0 & 0 \\
       0 & I_N\otimes K_1 & 0 \\
       0 & 0 & I_N\otimes K_2 \\
     \end{array}
   \right)x_c(t) \\
   =& K x_c(t).
\end{split}
\end{equation}
The if part of the proof is thus completed.

To show the only if part, we only need to note that, system (\ref{Eqlemc1}) is asymptotically stable only
if system (\ref{Thee2}) is asymptotically stable and  only if the digraph satisfies Assumption \ref{Ass2.6}.
\end{Proof}

\begin{Theorem}\label{Theo4.2}
Under Assumptions \ref{Ass2.3}-\ref{Ass2.5}, \ref{Ass2.2} and \ref{Ass2.7}, there exist matrices $K_1 \in \R^{m \times n_z}$, $K_2 \in \R^{m \times n}$ and $L \in \R^{n \times p}$ such that the cooperative robust output regulation problem is solved by the distributed dynamic output
feedback control law \eqref{ctr22} with $(G_1, G_2)$ being the minimal $p$-copy internal model of $S$ if and only if Assumption \ref{Ass2.6} is satisfied.
\end{Theorem}
\begin{Proof}
Let $x_c = \mbox{col}(x ,\bar{z},\bar{\xi} )$. Then the closed-loop system composed of the auxiliary  system \eqref{Csystemauxi} and the dynamic output feedback control law  \eqref{Ctr1n2}
is the same as the closed-loop system composed of the auxiliary augmented system \eqref{syscomp} and the dynamic output feedback control law \eqref{Ctr1n20} and
can be put into the following form:
\begin{equation}\label{cl30}
\begin{split}
x_c(t+1) &= \sum_{l=0}^{1} A_{cwl} x_c(t-\bar{r}_l) + B_{cw} v(t),\\
e_v(t) &=  C_{cw} x_c(t) + D_{cw} v(t), \\
\end{split}
\end{equation}
where $\bar{r}_0=0$, $\bar{r}_1=r$, and
\begin{equation*}\label{}
\begin{split}
A_{cw0} =& \left(
          \begin{array}{ccc}
            \bar{A} & 0 & 0\\
            \bar{G}_2 \bar{C}  & \bar{G}_1 & 0\\
            \bar{S}_2 \bar{C} & 0 & \bar{S}_1 \\
          \end{array}
        \right),\
A_{cw1} = \left(
          \begin{array}{cc}
             0 &
                    \begin{array}{cc}
                      \bar{B} \bar{K}_1 & \bar{B} \bar{K}_2 \\
                    \end{array}
                    \\
                \begin{array}{c}
                  0 \\
                  0 \\
                \end{array}    &
                                \begin{array}{c}
                                  0 \\
                                  \bar{S}_3 \\
                                \end{array}\\
          \end{array}
        \right), \\
B_{cw}= & \left(
            \begin{array}{c}
              \bar{E} \\
              \bar{G}_2 \bar{F}  \\
              \bar{S}_2 \bar{F}  \\
            \end{array}
          \right),
C_{cw} = \left(
             \begin{array}{ccc}
               \bar{C}  & 0 & 0\\
             \end{array}
           \right),
D_{cw}= \bar{F}. \\
\end{split}
\end{equation*}
Thus,  the nominal closed-loop system with $v$ set to $0$ is as follows:
\begin{equation}\label{eqqThee1}
\begin{aligned}
x_c(t+1) =& \left(
               \begin{array}{ccc}
                  I_N\otimes A  & 0 & 0 \\
                  H\otimes G_2 C  & I_N \otimes G_1 & 0 \\
                  H\otimes LC  &  0  & I_N\otimes A - H\otimes LC\\
               \end{array}
             \right) x_c(t) \\& + \left(
                              \begin{array}{ccc}
                                0&  I_N\otimes B  & H\otimes B   \\
                                0& 0 & 0\\
                                0 & I_N\otimes B & H\otimes B \\
                              \end{array}
                            \right) \left(
                                      \begin{array}{ccc}
                                        I_N \otimes K_2 & 0 & 0 \\
                                        0 & I_N \otimes K_1 & 0 \\
                                        0 & 0 & I_N \otimes K_2  \\
                                      \end{array}
                                    \right)
                            x_c(t-r),\\
\end{aligned}
\end{equation}
where $x_c=\mbox{col}(x,\bar{z},\bar{\xi})$ with $x=\mbox{col}(x_1,\dots,x_N)$, $\bar{z}=\mbox{col}(\bar{z}_1,\dots,\bar{z}_N)$ and $\bar{\xi}=\mbox{col}(\bar{\xi}_1,\dots,\bar{\xi}_N)$.

By Lemma \ref{Lem4.5}, there exist matrices $K_1 \in \R^{m \times n_z}$, $K_2 \in \R^{m \times n}$ and $L \in \R^{n \times p}$, such that system \eqref{eqqThee1} is asymptotically stable. The proof is thus completed by invoking Lemma \ref{lem3.1}.
\end{Proof}
\vspace{-6pt}
\section{EXAMPLE}\label{example}

Consider the discrete-time linear time-delay multi-agent systems of the form \eqref{system} with $N=4$, $r_{con}=1$, $\bar{A}_i=\left[\begin{array}{cc}
                    1&1+w_{i1}\\
                    0&1\\
                    \end{array}\right]$, $\bar{B}_i=\left[\begin{array}{cc}
                            1+w_{i2}\\
                            1
                            \end{array}\right]$, $\bar{E}_i=\left[\begin{array}{cc}
                            0&w_{i3}\\
                            0&i\\
                            \end{array}\right]$, $\bar{C}_i=\left[
     \begin{array}{cc}
       1 & 0 \\
     \end{array}
   \right]$, $F_i=\left[
     \begin{array}{cc}
       -1 & 0 \\
     \end{array}
   \right]$, $i=1,2,3,4$, and $v(t)$ is generated by the following exosystem:
\begin{equation}
v(t+1)=\left[\begin{array}{cc}
\cos{1}&\sin{1}\\
-\sin{1}&\cos{1}
\end{array}\right]v(t).
\end{equation}
The nominal system matrices are $A=\left[\begin{array}{cc}
1&1\\
0&1\\
\end{array}\right],$ $B=\left[\begin{array}{cc}1\\
1\\
\end{array}\right]$, $C=\left[\begin{array}{cc}
1&0\\
\end{array}\right],$ $F=\left[\begin{array}{cc}
-1&0\\
\end{array}\right],$ $E_i=\left[
    \begin{array}{cc}
      0&0\\
      0&i\\
    \end{array}
  \right],$ $i=1,2,3,4.$

\begin{figure}[H]
\begin{center}
\scalebox{0.4}{\includegraphics{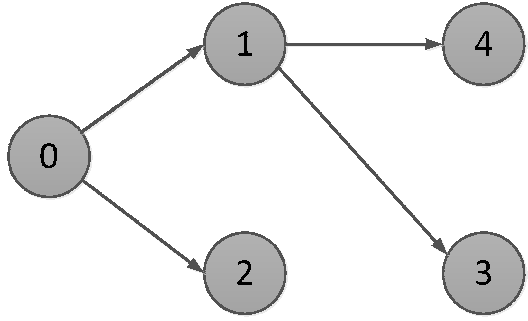}}\caption{The network topology}\label{graph}
\end{center}
\end{figure}

The communication network topology is described in Fig. \ref{graph}. The matrix $H$ associated with the digraph $\bar{\mathcal{G}}$ is \begin{equation*}
\begin{split}
H=\left[ \begin{array}{cccc}
            1&0&0&0\\
            0&1&0&0\\
            -1&0&1&0\\
            -1&0&0&1\\
            \end{array}\right],
\end{split}
\end{equation*} whose eigenvalues  are $\{1,1,1,1\}$.

It is easy to verify that Assumptions \ref{Ass2.3}-\ref{Ass2.7} are satisfied. Therefore, by Theorems \ref{Theo4.1} and \ref{Theo4.2}, the cooperative robust output regulation problem for this example can be solved by the distributed control laws of the form \eqref{ctr1tr} and \eqref{ctr22}.

\noindent (1) Distributed dynamic state feedback control:

The distributed dynamic state feedback control law is given as in \eqref{ctr1tr} with $i=1,\dots,4$, and

\begin{equation}\label{gg1}
G_1=\left[
\begin{array}{cc}
 \cos 1 & \sin 1  \\
  -\sin 1 & \cos 1 \\
  \end{array}
  \right]\ \text{and}\
G_2=\left[
             \begin{array}{cc}
               0 \\
               1 \\
             \end{array}
           \right].
\end{equation}
Assume the communication delay $r_{com}=1$.

Denote $ {A}_c=\left(
                      \begin{array}{cc}
                        A & 0 \\
                        G_2C & G_1 \\
                      \end{array}
                    \right)$ and $ {B}_c=\left(
                                           \begin{array}{c}
                                             B \\
                                             0 \\
                                           \end{array}
                                         \right)
                    $. By Lemma \ref{Lem4.3ap}, the desirable feedback gain is 
\begin{equation}\label{gg2}
K=(K_x,K_z)=-\nu_1^{-1} R^{-1}B_c^TPA_c^{r+1},
\end{equation}
where $R=I_m+B_c^TPB_c$, $r=2$ and $P$ is the positive definite solution of the parametric DARE
\begin{equation}\label{}
A_c^T PA_c-P  -A_c^TPB_cR^{-1}B_c^T PA_c=-\gamma P,
\end{equation}
where $\gamma$ is some sufficiently small positive number. Then, $$K=\left[\begin{array}{cccc} 0.1292&-0.1788&-0.0659&-0.1597 \end{array}\right].$$

\begin{figure}[H]
\begin{center}
\scalebox{0.8}{\includegraphics{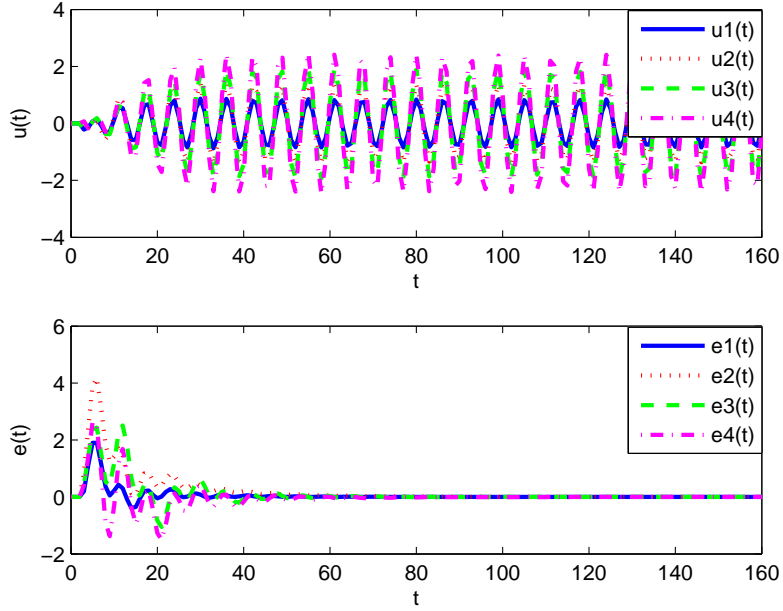}}\caption{Control inputs and tracking errors under distributed dynamic state feedback control}\label{inputst}
\end{center}
\end{figure}

With random initial conditions, Fig. \ref{inputst} shows the control inputs of the system which are bounded, and the tracking errors of the followers which tend to zero asymptotically. The system uncertainties are $w\!=\!(0.1, 0.2, 0.3, 0.4, 0.1, 0.2, 0.3, 0.4, 0.5, 0.6, 0.7, 0.8)^T\!$, $\nu_1=1$ and $\gamma=0.08$.

\noindent (2) Distributed dynamic output feedback control:

The distributed dynamic output feedback control law is given as in \eqref{ctr22} with $i=1,\dots,4$, and $(G_1,G_2)$, $(K_1,K_2)=(K_z,K_x)$ defined in \eqref{gg1} and \eqref{gg2}, respectively. By Lemma \ref{Lem4.5}, the desirable observer gain is $L=-K_l^T$, where
\begin{equation}
K_l=-\nu_2^{-1}R_l^{-1}CP_l(A^T)^{r+1}
\end{equation}
with $R_l=I_m+CP_lC^T,$ and $P_l$ is the positive-definite solution of the parametric DARE
\begin{equation}
AP_lA^T-P_l-AP_lC^TR_l^{-1}CP_lA^T=-\gamma_lP_l.
\end{equation}
Then $L=\mbox{col}(0.72, 0.0648)$.

\begin{figure}[H]
\begin{center}
\scalebox{0.8}{\includegraphics{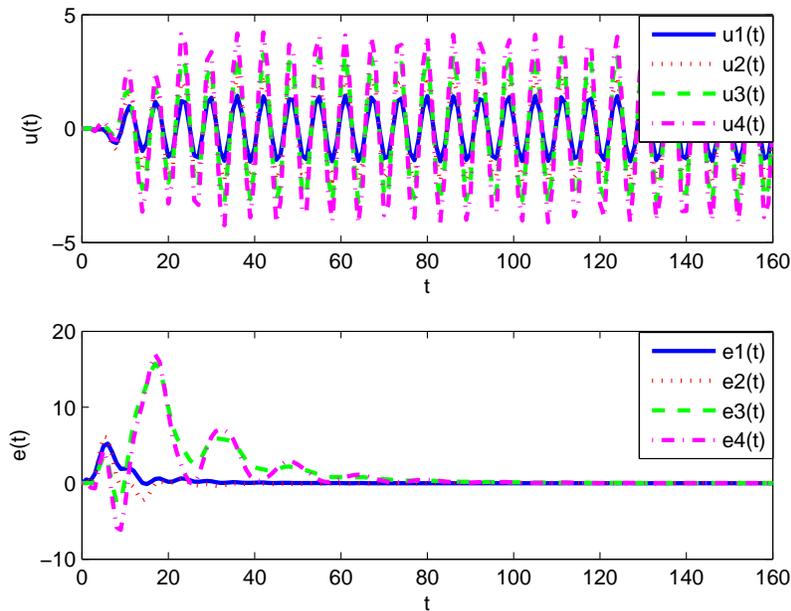}}\caption{Control inputs and tracking errors under distributed dynamic output feedback control}\label{inputop}
\end{center}
\end{figure}

With random initial conditions, Fig. \ref{inputop} shows the control inputs of the system which are bounded, and the tracking errors of the followers which tend to zero asymptotically. The system uncertainties are $w\!=\!(0.1,0.2,0.3,0.4,0.1,0.2,0.3,0.4,0.5,0.6,0.7,0.8)^T\!$, $\nu_2=0.5$ and $\gamma_l=0.18$. 

\section{CONCLUSION}\label{concl}

 This paper has investigated the cooperative robust output regulation problem for discrete-time linear multi-agent systems with input and communication delays by dynamic distributed internal model approach. We have established the solvability conditions for the problem via both the distributed state feedback control and the distributed output feedback control. Future work will focus on stabilization of auxiliary  augmented system by other approaches with multiple input and communication delays.


\section*{APPENDIX}
\vspace{-2pt}
\begin{Lemma}(Lemma 10.2 in \cite{Bzhoubook})\label{Lem4.1ap}
Assume that $(A,B)$ is controllable, all the eigenvalues of $A$ are on the unit circle, and $\lambda \in \mathcal{C}$
with $\mbox{Re}(\lambda) > 0$. Let $K= -\mu^{-1} R^{-1}B^T P(\gamma)A^{r+1}$ where $R=I_m+B^TP(\gamma)B,$ and $r \in \mathbb{Z}^+$, and $ P(\gamma)$ is the unique positive definite solution to the
parametric discrete-time algebraic Riccati equation (DARE)
\begin{equation}\label{DARE}
 {A}^T P(\gamma)A - P(\gamma) - A^TP(\gamma) {B}{R}^{-1} {B}^TP(\gamma)A=-\gamma P(\gamma).
\end{equation}
Then, for any $0< \mu \leq \mbox{Re}(\lambda)$, there exists a positive scalar $ \gamma^*$ such that the system
\begin{equation}
x(t+1)=Ax(t) +\lambda B K x(t-r)
\end{equation}
is asymptotically stable for all $\gamma \in (0, \gamma^*]$.
\end{Lemma}

\begin{Remark}\label{Rem4.3ap}
As pointed out in \cite{Bzhoubook}, the assumption in Lemma \ref{Lem4.1ap} that all the eigenvalues of the matrix $A$ have modulus $1$ can be relaxed to the assumption that all the eigenvalues of $A$ have modulus equal to or smaller than $1$, and the assumption that $(A, B)$ is controllable can be relaxed to $(A, B)$ is stablizable.
\end{Remark}

\begin{Lemma} \label{Lem4.3ap}
 Consider the system of the form
\begin{equation}\label{syslem1}
x_i(t+1) = A x_i(t) + \lambda_i B u_i(t-r),\ i=1,\dots,N,
\end{equation}
where $x_i \in \R^{n}$, $u_i \in \R^{m}$, ${A} \in \R^{n \times n}$,  ${B} \in \R^{n \times m}$, and $\lambda_i \in \mathcal {C}$ with $Re\{\lambda_i\} >0$.
Suppose all the  eigenvalues of $A$ have modulus equal to or smaller than $1$, and $(A, B)$ is stablizable. Then,
 there exists a matrix $K \in \R^{m \times n}$ such that the state feedback control law $u_i(t)=Kx_i(t)$,
$i=1,\dots,N$, asymptotically stabilize all subsystems of the system \eqref{syslem1}.
\end{Lemma}

\begin{Proof}
Since  ${A}$ has no eigenvalues with modulus greater than $1$,  there exists a non-singular matrix ${T}$ such that
\begin{equation*}
{T}A {T}^{-1}=\left(\begin{array}{cc}
A_1&0\\
0&A_2\\
\end{array}\right),\ \ {T}B=\left(\begin{array}{cc}
B_{1}\\
 B_{2} \\
\end{array}\right),
\end{equation*}
where all the eigenvalues of $A_1 \in \R^{n_1 \times n_1} $ have modulus $1$ and all the eigenvalues of $A_2$ have modulus smaller than $1$. Moreover, $\left(A_1, B_1\right)$ is controllable. Let $\hat{x}_i = {T} x_i = \mbox{col } (\hat{x}_{i1}, \hat{x}_{i2})$ with $\hat{x}_{i1} \in \R^{n_1} $. Then \eqref{syslem1} is transformed to the following:
\begin{equation}
\begin{split}
\hat{x}_{i1}(t+1) =&A_1 \hat{x}_{i1} (t) +   \lambda_i  B_{1} u_i(t-r), \\
\hat{x}_{i2}(t+1) =&A_2 \hat{x}_{i2} (t) +   \lambda_i  B_{2} u_i(t-r).
\end{split}
\end{equation}
By Lemma \ref{Lem4.1ap}, there exists a $\gamma^*>0$ such that, for any $\gamma \in (0,\gamma^*)$,
the following parametric DARE,
\begin{equation}
A_1^TP_1A_1-P_1-A_1^TP_1B_1 R_1^{-1}B_1^TP_1A_1=-\gamma P_1,
\end{equation}
where $R_1 = I_{m_1}+B_1^TP_1B_1$ has a unique positive-definite solution $P_1$. Moreover, let $K_1 = - \mu_1^{-1}R_1^{-1}B_1^TP_1A_1^{r+1}$, where $0<\mu_1\leq Re\{\lambda_i\}$. Then,
\begin{equation}\label{tpfcon2}
\hat{x}_{i1} (t+1) =A_1 \hat{x}_{i1} (t)  + \lambda_i B_{1} K_1  \hat{x}_{i1}  (t-r) \\
\end{equation}
is asymptotically stable. Since $A_2$ is Schur, under the control $u_i (t) = K_1  \hat{x}_{i1}  (t)$, $\hat{x}_{i2}(t)$ also tends to zero as $t$ tends to infinity. Thus,   the state feedback control law $u_i(t)=Kx_i(t)$ where ${K}=\big(K_1, 0\big) {T}$,
$i=1,\dots,N$, asymptotically stabilize all subsystems of the system \eqref{syslem1}.
\end{Proof}


%
%
%


\begin{thebibliography}{90}

%

\bibitem{Da}
Davison EJ. The robust control of a servomechanism problem for linear time-invariant multivariable systems. \emph{IEEE Transactions on Automatic Control}, 1976; \textbf{21}(1), 25--34. 

\bibitem{Fr}
Francis BA. The linear multivariable regulator problem. \emph{SIAM Journal on Control and Optimization}, 1977; \textbf{15}(3): 486--505.

\bibitem{FW}
Francis BA, Wonham WM. The internal model principle of control theory. \emph{Automatica}, 1976; \textbf{12}(5): 457--465.

\bibitem{Su1}
Su Y, Huang J. Cooperative output regulation of linear multi-agent systems. \emph{IEEE Transactions on Automatic Control}, 2012; \textbf{57}(4): 1062--1066.

\bibitem{Su3}
Su Y, Hong Y, Huang J. A general result on the robust cooperative output regulation for linear uncertain multi-agent systems. \emph{IEEE Transactions on Automatic Control}, 2013; \textbf{58}(5): 1275--1279.

\bibitem{Jiang1}
Wang X, Hong Y, Huang J, Jiang Z. A distributed control approach to a robust output regulation problem for multi-agent linear systems. \emph{IEEE Transactions on Automatic Control}, 2011; \textbf{55}(12): 2891--2895.

\bibitem{CCC14}
Lu M, Huang J. Cooperative output regulation problem for linear time-delay multi-agent systems under switching network. \emph{ Proc. 33rd Chinese Control Conference}, pp. 3515--3520, Nanjing, China, 2014.

\bibitem{Lu2015}
Lu M, Huang J. Internal model approach to cooperative robust output regulation for linear uncertain time-delay multi-agent systems. \emph{Internaltional Journal of Robust and Nonlinear Control}, submitted, available in  https://arxiv.org/abs/1508.04207 

\bibitem{IET}
Yan Y, Huang J. Cooperative output regulation of discrete-time linear time-delay multi-agent systems. \emph{IET Control Theory and Applications}, 2016; \textbf{10}(16): 2019--2026.

\bibitem{JSSC}
Yan Y, Huang J. Robust output regulation problem for discrete-time linear systems with both input and communication delays. \emph{Journal of Systems Science and Complexity}, 2017; \textbf{30}(1): 68--85.

\bibitem{Liang2014}
Liang H, Zhang H, Wang Z, Wang J. Consensus robust output regulation of discrete-time linear multi-agent systems. \emph{IEEE/CAA Journal of Automatica Sinica}, 2014; \textbf{1}(2): 204--209.

\bibitem{chen1}
Chen J, Sun J, Liu G, Rees D.  New delay-dependent stability criteria for neural networks with time-varying interval delay. \emph{Physics Letters A}, 2010; \textbf{374}(43): 4397--4405.

\bibitem{chen2}
Sun J, Chen J, Liu G, Rees D. On robust stability of uncertain neutral systems with discrete and distributed delays. \emph{2009 American Control Conference}, pp. 5469--5473, St. Louis, Missouri, USA, 2009.

\bibitem{you2013}
Hengster-Movrica K, You K, Lewis F, Xie L. Synchronization of discrete-time multi-agent systems on graphs using
Riccati design. \emph{Automatica}, 2013; \textbf{49}(2): 414--423.


\bibitem{Graph}
Godsil C, Royle G.(2001). \emph{Algebraic Graph Theory}. New York: Springer-Verlag.







\bibitem{J.H}
Huang, J.(2004). \emph{Nonlinear Output Regulation: Theory and Applications}, Philadelphia, USA: SIAM.

\bibitem{Bzhoubook}
Zhou, B.(2014). \emph{Truncated Predictor Feedback for Time-delay Systems}, Springer-Verlag Berlin Heidelberg.




%
%























\end{thebibliography}
\end{document}